%% file: blocsarxiv.tex
\documentclass[11pt,twoside]{article}
\usepackage[T1]{fontenc}
\usepackage{amsmath,amsfonts,amssymb,amsthm,fullpage,bbm}
\usepackage{fancyhdr,graphicx,color}
\usepackage{epsfig}

\newtheorem{thm}{Theorem}
\newtheorem{lem}{Lemma}

\theoremstyle{definition}

\theoremstyle{remark}

\title{\huge On the small maximal flows in first passage percolation} 
\author{\Large Marie THERET}
\date{}

\begin{document}
\maketitle

We consider the standard first passage percolation on $\mathbb{Z}^{d}$: with each edge of the lattice we associate a random capacity. We are interested in the maximal flow through a cylinder in this graph. Under some assumptions Kesten proved in 1987 a law of large numbers for the rescaled flow. Chayes and Chayes established that the large deviations far away below its typical value are of surface order, at least for the Bernoulli percolation and cylinders of certain height. Thanks to another approach we extend here their result to higher cylinders, and we transport this result to the model of first passage percolation.

\section{Definitions and main result}

We use the notations introduced in \cite{Kesten:StFlour} and \cite{Kesten:flows}. Let $d\geq2$. We consider the graph $(\mathbb{Z}^{d}, \mathbb E ^{d})$ having for vertices $\mathbb Z ^{d}$ and for edges $\mathbb E ^{d}$ the set of the pairs of nearest neighbors for the standard $L^{1}$ norm. With each edge $e$ in $\mathbb{E}^{d}$ we associate a random variable $t(e)$ with values in $\mathbb{R}^{+}$. We suppose that the family $(t(e), e \in \mathbb{E}^{d})$ is independent and identically distributed, with a common distribution function $F$. More formally, we take the product measure $\mathbb {P}$ on $\Omega= \prod_{e\in \mathbb{E}^{d}} [0, \infty[$, and we write its expectation $\mathbb{E}$. We interpret $t(e)$ as the capacity of the edge $e$; it means that $t(e)$ is the maximal amount of fluid that can go through the edge $e$ per unit of time. For a given realization $(t(e),e\in \mathbb{E}^{d})$ we denote by $\phi_{\vec{k},m} = \phi_{B}$ the maximal flow through the box
$$ B(\vec{k},m) \, = \, \prod _{i=1}^{d-1} [0,k_{i}] \times [0,m] \, ,$$
where $\vec{k}=(k_{1},...,k_{d-1}) \in \mathbb{Z}^{d-1}$, from its bottom
$$F_{0}\, = \, \prod _{i=1}^{d-1} [0,k_{i}] \times \{0\} $$
to its top
$$F_{m}\, = \, \prod _{i=1}^{d-1} [0,k_{i}] \times \{m\} \, .$$
Let us define this quantity properly. We remember that $\mathbb{E}^{d}$ is the set of the edges of the graph. An edge $e\in \mathbb{E}^{d}$ can be written $e=\langle x,y \rangle$, where $x$, $y\in \mathbb{Z}^{d}$ are the endpoints of $e$. We will say that $e=\langle x,y\rangle$ belongs to a subset $A$ of $\mathbb{R}^{d}$ ($e\in A$) if the segment joining $x$ to $y$ (eventually excluding these points) is included in $A$. Now we define $\widetilde{\mathbb{E}}^{d}$ as the set of all the oriented edges, i.e. an element $\widetilde{e}$ in $\widetilde{\mathbb{E}}^{d}$ is an ordered pair of vertices. We denote an element $\widetilde{e} \in \widetilde{\mathbb{E}}^{d}$ by $\langle \langle x,y \rangle \rangle$, where $x$, $y \in \mathbb{Z}^{d}$ are the endpoints of $\widetilde{e}$ and the edge is oriented from $x$ towards $y$. We consider now the set $\mathcal{S}$ of all pairs of functions $(g,o)$, with $g:\mathbb{E}^{d} \rightarrow \mathbb{R}^{+}$ and $o:\mathbb{E}^{d} \rightarrow \widetilde{\mathbb{E}}^{d}$ such that $o(\langle x,y\rangle ) \in \{ \langle \langle x,y\rangle \rangle , \langle \langle y,x \rangle \rangle \}$, satisfying
\begin{itemize}
\item for each edge $e$ in $B$ we have
$$0 \,\leq\, g(e) \,\leq\, t(e) \,,$$
\item for each vertex $v$ in $B \smallsetminus F_{m}$ we have
$$ \sum_{e\in B\,:\, o(e)=\langle\langle v,\cdot \rangle \rangle} g(e) \,=\, \sum_{e\in B\,:\, o(e)=\langle\langle \cdot ,v \rangle \rangle} g(e) \,. $$ 
\end{itemize}
A couple $(g,o) \in \mathcal{S}$ is a possible stream in $B$: $g(e)$ is the amount of fluid that goes through the edge $e$, and $o(e)$ gives the direction in which the fluid goes through $e$. The two conditions on $(g,o)$ express only the fact that the amount of fluid that can go through an edge is bounded by its capacity, and that there is no loss of fluid in the cylinder. With each possible stream we associate the corresponding flow
$$ flow (g,o) \,=\, \sum_{ u \notin F_{m} \,,\,  v \in F_{m} \,:\, \langle
  u,v\rangle \in \mathbb{E}^{d}\cap B} g(\langle u,v\rangle) \mathbb{I}_{o(\langle u,v\rangle) = \langle\langle u,v \rangle\rangle} - g(\langle u,v\rangle) \mathbb{I}_{o(\langle u,v\rangle) = \langle\langle v,u \rangle\rangle} \,. $$
This is the amount of fluid that crosses the cylinder $B$ if the fluid respects the stream $(g,o)$. The maximal flow through the cylinder $B$ is the supremum of this quantity over all possible choices of stream
$$ \phi_{B} \,=\, \phi_{\vec{k},m} \,=\, \sup_{(g,o) \in \mathcal{S}} \, flow (g,o) \,.$$

We note $p_{c}(d)$ the critical value of the parameter of the Bernoulli percolation in dimension $d$. We will prove the following result:

\begin{thm}
\label{casgen}
We suppose that
$$ F(0) \, < \, 1-p_{c}(d) \,.$$
There exist a positive constant $\varepsilon_{0}$, depending only on $d$ and $F$, and a positive constant $C$, depending only on $d$, such that for any function $h: \mathbb{N} \rightarrow \mathbb{N}$ satisfying
$$ \lim_{n\rightarrow \infty} \frac{ \ln h(n)}{n^{d-1}} \, = \, 0  $$
we have
$$\forall \varepsilon < \varepsilon_{0} \qquad \liminf_{n\rightarrow \infty} -\frac{1}{n^{d-1}} \ln \mathbb{P} \left[ \phi_{(n,...,n),h(n)} \leq \varepsilon n^{d-1} \right] \, \geq C \, > \, 0 \, .  $$
\end{thm}

The condition $F(0)<1-p_{c}$ is necessary for this result to hold. Indeed, Yu Zhang (see \cite{Zhang}) proved in dimension $3$ that for a function $F$ satisfying
$$ F(0) \,=\, 1-p_{c} \qquad and\qquad \int_{[0,+\infty [} x dF(x)\,<\,\infty$$
we have
$$ \lim_{k,l,m\rightarrow \infty} \frac{\phi_{(k,l),m}}{kl} \,=\,0 \,. $$

The spirit of this result is not new, Chayes and Chayes proved in \cite{Chayes} (see Lemma $3.3$) the following theorem:
\begin{thm}
\label{thmChayes}
We suppose that the capacity $t$ of each edge follows a Bernoulli law of parameter $p$ satisfying $p>p_{c}$. Then there exist positive constants $\widetilde{\varepsilon}$, $\widetilde{C}$ such that
$$ \mathbb{P} \left[ \phi_{(n,...,n),n} \geq \widetilde{\varepsilon} n^{d-1} \right] \, \geq \, 1-e^{-\widetilde{C} n^{d-1}}  $$
for $n$ sufficiently large. 
\end{thm}
To prove it they divide the cylinder into thin layers, compare each one of them to objects of dimension $2$ and use the results of \cite{Aizenman-Chayes}. Because of the passage in dimension $2$, it seems to us that this proof can only be extended to cylinders $B((n,...,n),h(n))$ with a height satisfying $\lim_{n\rightarrow \infty} \ln h(n) / n =0$. This is the constraint we have in dimension $2$, but not in higher dimensions. Actually the condition $\lim_{n\rightarrow \infty} \ln h(n) / n^{d-1} =0$ is the good one, in the sense that in the model of Bernoulli percolation if $h(n) = \exp (k n^{d-1})$ for a constant $k$ sufficiently large, the maximal flow $\phi_{(n,...,n),h(n)}$ tends to $0$ almost surely. Indeed if the $n^{d-1}$ vertical edges of the cylinder that intersect one fixed horizontal plane have all $0$ for capacity then $\phi_{(n,...,n),h(n)}=0$. By independence and translation invariance of the model, we obtain, for $k$ large enough
$$ \mathbb{P} \left[ \phi_{(n,...,n),h(n)} \neq 0 \right] \,\leq \, \left[ 1-(1-p)^{n^{d-1}} \right]^{h(n)} \, \rightarrow_{n\rightarrow \infty} \, 0 \,. $$

The proof of Theorem \ref{casgen} is based on the coarse graining techniques of Pisztora (see \cite{Pisztora}). Actually we don't need estimates as strong as those of Pisztora for the renormalization scheme. We will use a weaker version of these results as in \cite{Cerf:StFlour}. Moreover we won't use the general stochastic domination inequality (see \cite{Pisztora}, \cite{Liggett}), it is sufficient here to use a partition of the space into equivalence classes to get rid of problems of dependence between random variables, as we will see in section \ref{equiv}.

We will first study two particular cases of this result in the model of Bernoulli percolation, that will allow us to deal very simply with the proof of the main theorem in general first passage percolation.


\section{max-flow min-cut theorem}

The definition of the flow is not easy to deal with. The maximal flow $\phi_{B}$ can be expressed differently thanks to the max-flow min-cut theorem (see \cite{Bollobas}). We need some definitions.

A path on a graph ($\mathbb{Z}^{d}$ for example) from $v_{0}$ to $v_{n}$ is a sequence $(v_{0}, e_{1}, v_{1},..., e_{n}, v_{n})$ of vertices $v_{0},..., v_{n}$ alternating with edges $e_{1},..., e_{n}$ such that $v_{i-1}$ and $v_{i}$ are neighbors in the graph, joined by the edge $e_{i}$, for $i$ in $\{1,..., n\}$. Two paths are said disjoint if they have no common edge.

A set $E$ of edges of $B(\vec{k},m)$ is said to separate $F_{0}$ from $F_{m}$ in $B(\vec{k},m)$ if there is no path from $F_{0}$ to $F_{m}$ in $B(\vec{k},m) \smallsetminus E$. We call $E$ an $(F_{0},F_{m})$-cut if $E$ separates $F_{0}$ from $F_{m}$ in $B(\vec{k},m)$ and if no proper subset of $E$ does. With each set $E$ of edges we associate the variable
$$ V(E)\, = \, \sum_{e\in E} t(e) \, .$$
The max-flow min-cut theorem states that
$$ \phi_{B} \, = \, \min \{ \, V(E) \, | \, E \,\, is \,\, an \,\, (F_{0},F_{m})-cut \, \} \, .$$

In the special case where $t(e)$ belongs to $\{0,1\}$, i.e. the law of $t$ is a Bernoulli law, the flow has an other simple expression. In this case, let us consider the graph obtained from the initial graph $\mathbb{Z}^{d}$ by removing all the edges $e$ with $t(e)=0$. Menger's theorem (see \cite{Bollobas}) states that the minimal number of edges in $B(\vec{k},m)$ that have to be removed from this graph to disconnect $F_{0}$ from $F_{m}$ is exactly the maximal number of disjoint paths that connect $F_{0}$ to $F_{m}$. By the max-flow min-cut theorem, it follows immediately that the maximal flow in the initial graph through $B$ from $F_{0}$ to $F_{m}$ is exactly the maximal number of disjoint open paths from $F_{0}$ to $F_{m}$, where a path is open if and only if the capacity of all its edges is one.


\section{Bernoulli percolation for a parameter $p$ near $1$}

We consider that the capacity $t$ of each edge follows the Bernoulli law of parameter $p$, with $p = \mathbb{P}[t=1]$ as close to $1$ as we will need. Remember that here the maximal flow through a cylinder $B$ is the maximal number of disjoint open paths from the bottom to the top of $B$. We will first prove the following theorem:

\begin{thm}
\label{bernoulli1}
For all $\varepsilon$ in $[0,1[$, there exist $p_{0}(\varepsilon,d)<1$ and a constant $C'$ depending only on the dimension $d$ such that for any function $h: \mathbb{N} \rightarrow \mathbb{N}$ satisfying
$$ \lim_{n\rightarrow \infty} \frac{ \ln h(n)}{n^{d-1}} \, = \, 0  $$
and for all $p \geq p_{0}$ we have
$$ \liminf_{n\rightarrow \infty} -\frac{1}{n^{d-1}} \ln \mathbb{P} \left[ \phi_{(n,...,n),h(n)} \leq \varepsilon n^{d-1} \right] \, \geq C' \, > \, 0 \, .  $$
\end{thm}

To simplify the notations during the proof of this theorem, we define
$$ \alpha(\varepsilon) \,=\, \mathbb{P} \left[ \phi_{(n,...,n),h(n)} \leq \varepsilon n^{d-1}  \right] \,.$$
Thanks to the max-flow min-cut theorem, we know that
$$ \alpha(\varepsilon) \,=\, \mathbb{P} \left[\, there \,\, exists \,\, a \,\, (F_{0},F_{h(n)})-cut \,\, E \,\, satisfying \,\, V(E)\leq \varepsilon n^{d-1}\,  \right] \,.$$
We need to define a notion of $\diamond$-connection. We associate with each
edge $e$ a plaquette $\mathcal{P}(e)$ which is the only unit square of
$\mathbb{R}^{d}$ of the form $\mathcal{P}_{i} + (n_{1},...,n_{d})$ that
intersects $e$ in its middle, where $(n_{1},...,n_{d}) \in \mathbb{Z}^{d}$
and $\mathcal{P}_{i} = [-1/2,1/2]^{i-1} \times \{1/2\} \times
[-1/2,1/2]^{d-i}$ for $1\leq i \leq d$. We say that two edges $e_{1}$ and
$e_{2}$ are $\diamond$-connected if and only if $\mathcal{P}(e_{1}) \cap
\mathcal{P}(e_{2}) \neq \emptyset$. According to Kesten (see
\cite{Kesten:flows}) a $(F_{0},F_{h(n)})$-cut is
$\diamond$-connected. Moreover, it is obvious that a cut contains at least
$n^{d-1}$ edges (to cut the $n^{d-1}$ possible vertical paths). In
particular, if we consider a fixed vertical path and if we denote by
$(e_{i},i=1,...,h(n))$ the edges of this path, a $(F_{0},F_{h(n)})$-cut $E$
must contain at least one of these $e_{i}$. We can then find a subset $E'$
of $E$ which contains exactly $n^{d-1}$ edges, including one of these
$e_{i}$, and which is $\diamond$-connected. Obviously if $V(E) \leq
\varepsilon n^{d-1}$ then $V(E') \leq \varepsilon n^{d-1}$. Finally we can
relax the constraint for $E'$ to be in $B$, and by translation invariance
of the model we can suppose that $E'$ contains a determined edge $e_{0}$. We deduce from these remarks that
\begin{align*}
\alpha(\varepsilon) & \,\leq\, \sum_{i=1}^{h(n)} \, \mathbb{P} \left[ \,there \,\, exists \, a \,\, (F_{0},F_{h(n)})-cut \,\, E \,\, such \,\, that \,\, V(E)\leq \varepsilon n^{d-1} \,\, and \,\, e_{i}\in E \, \right] \\
& \,\leq\, h(n)\, \mathbb{P} \left[\, \begin{array}{c} there \,\, exists \,\, a \,\, \diamond -connected \,\, set \,\, E' \,\, of \,\, n^{d-1} \,\, edges\\  such \,\, that\,\, V(E')\leq \varepsilon n^{d-1} \,\, and \,\, e_{0}\in E'\end{array}  \, \right] \\
& \,\leq\, h(n) \sum_{A} \mathbb{P} \left[ \sum_{e\in A} t(e) \leq \varepsilon n^{d-1} \right] \, , \\
\end{align*}
where the sum is over the $\diamond$-connected sets $A$ of $n^{d-1}$ edges including $e_{0}$. We know (see \cite{Kesten:flows}) that there exists a constant $c>1$ depending only on the dimension $d$ such that the number of such possible sets $A$ is bounded by $c^{n^{d-1}}$. We deduce then, thanks to the exponential Chebyshev inequality, that
\begin{align*}
\forall \lambda >0 \qquad \alpha(\varepsilon) & \,\leq\, h(n) c^{n^{d-1}}e^{\lambda\varepsilon n^{d-1}} \mathbb{E} \left[ e^{-\lambda t} \right] ^{n^{d-1}} \\
& \, \leq \, \exp \left(-n^{d-1} \left[ -\frac{\ln h(n)}{n^{d-1}} - \ln c + \lambda ( 1-\varepsilon) - \ln (p+(1-p) e^{\lambda}) \right] \right) \,.\\
\end{align*}
We choose $\lambda$ such that
$$ \lambda (1-\varepsilon) \, \geq \, 3 \ln c \,,$$
and then $p_{0} <1$ (depending on $d$ and $\varepsilon$) such that for all $p\geq p_{0}$ we have
$$ \ln \left( p+(1-p) e^{\lambda} \right) \leq \ln c \,.$$
We conclude that
$$ \forall p\geq p_{0} \qquad \alpha(\epsilon) \, \leq \, \exp \left(-n^{d-1} \left[ \ln c - \frac{\ln h(n)}{n^{d-1}}  \right] \right) \,.$$
This ends the proof of theorem \ref{bernoulli1}.


\section{Bernoulli percolation}

We consider now that the law of $t$ is a Bernoulli law with a fixed parameter $p>p_{c}$. We will prove the following result:
\begin{thm}
\label{bernoulli}
For any $p>p_{c}$, there exist a positive $\varepsilon_{0}$ (depending on $d$ and $p$) and a positive constant $C''$ (depending only on the dimension $d$) such that for any function $h: \mathbb{N} \rightarrow \mathbb{N}$ satisfying
$$ \lim_{n\rightarrow \infty} \frac{ \ln h(n)}{n^{d-1}} \, = \, 0  $$
we have
$$ \forall \varepsilon \leq \varepsilon_{0} \qquad \liminf_{n\rightarrow \infty} -\frac{1}{n^{d-1}} \ln \mathbb{P} \left[ \phi_{(n,...,n),h(n)} \leq \varepsilon n^{d-1} \right] \, \geq C'' \, > \, 0 \, .  $$
\end{thm}

The proof of this theorem is based on the coarse graining techniques of Pisztora (see \cite{Pisztora}, \cite{Cerf:StFlour}). The idea is to use a renormalization scheme: instead of looking at what happens for each edge, we try to understand what are the typical properties of the edges in a box, and to deduce some properties for the entire graph.


\subsection{Coarse graining}

Let $\Lambda$ be a box. We define its inner vertex boundary as
$$ \partial^{in} \Lambda \,=\, \{ x \in \Lambda \,|\, \exists \, y \notin \Lambda \,, \,\, |x-y|=1 \} \,. $$
An open cluster within $\Lambda$ is said crossing for $\Lambda$ if it
intersects each of the $2d$ faces of $\partial^{in}\Lambda$. The diameter
of a set $A$ is given by $diam(A)=\max_{i=1...d} \sup_{x,y\in A}
|x_{i}-y_{i}|$. We now consider the event
$$ U(\Lambda) \,=\, \{\,there\,\,exists\,\,an\,\,open\,\,crossing\,\,cluster\,\,in\,\,\Lambda  \,\}  $$
and, for $m$ less than or equal to the diameter of $\Lambda$,
$$ W(\Lambda,m)\,=\, \{\,there \,\,exists\,\, a\,\,unique\,\,open\,\,
cluster \,\,in\,\,\Lambda \,\, with \,\, diameter\,\,\geq m \, \} \,.  $$
Let $\Lambda(n)$ be the square box $]-n/2,n/2]^{d}$. We know that
\begin{lem}
\label{evtR}
For all dimension $d\geq 2$ and for all $p> p_{c}$, we have
$$ \lim_{n\rightarrow \infty} \mathbb{P} [U(\Lambda(n))] \,=\,1 \,.$$
Moreover, there exists a finite constant $\gamma$ (depending on $d$ and $p$) such that
$$ \lim_{n\rightarrow \infty} \mathbb{P} [W(\Lambda(n),\gamma\ln n)] \,=\,1 \,.$$
\end{lem}
For a proof, see \cite{Cerf:StFlour}.
In particular, this lemma implies that
$$ \lim_{n \rightarrow \infty} \mathbb{P} ( W(\Lambda (n),n/3)) \,=\,1
\,.$$

To use this estimate, we will rescale the lattice. Let $K$ be a positive integer. We divide $\mathbb{Z}^{d}$ into small boxes called blocks of size $K$ in the following way. For $\underline{x}=(\underline{x}_{1},...,\underline{x}_{d}) \in \mathbb{Z}^{d}$, we define the block indexed by $\underline{x}$ as
$$ B_{K}(\underline{x}) \,=\, K \underline{x} + \Lambda(K) \,,  $$
where $K\underline{x}$ is the vertex $(K\underline{x}_{1},...,K\underline{x}_{d})$. We remark that the blocks partition $\mathbb{R}^{d}$. Let $A$ be a region of $\mathbb{R}^{d}$, we define the rescaled region $\underline{A}_{K}$ as
$$ \underline{A}_{K} \,=\, \left\{ \underline{x} \in \mathbb{Z}^{d} \,|\, B_{K}(\underline{x}) \cap A \neq \emptyset \right\} \,.$$
For $\underline{x}\in \mathbb{Z}^{d}$, we define next a neighborhood of the block $B_{K}(\underline{x})$, called the event-block, as
$$ B'_{K}(\underline{x}) \,=\, \bigcup_{\underline{u} } B_{K}(\underline{u}) \,,$$
where the union is over the vertices
$\underline{u}=(\underline{u}_{1},...,\underline{u}_{d})\in \mathbb{Z}^{d}$
satisfying $\max_{1\leq i \leq d} |\underline{x}_{i} - \underline{u}_{i}|
\leq 1$. Finally we define the block process $(X_{K}(\underline{x}), \underline{x}\in\mathbb{Z}^{d})$ as
$$ \forall \underline{x}\in\mathbb{Z}^{d} \qquad X_{K}(\underline{x}) \,=\,
\mathbbm{1}_{U(B_{K}(\underline{x}))} \times
\mathbbm{1}_{W(B_{K}(\underline{x}), \frac{K}{3})} \times \prod_{\underline{y} \in Y}
\mathbbm{1}_{W(B_{K}(\underline{x}) + \underline{y}, \frac{K}{3})} \,,$$
where $Y= \{ (\pm K/2,0,...,0),(0,\pm K/2,0,...,0),...,(0,...,0,\pm K/2) \}$.
We say that the event-block $B'_{K}(\underline{x})$ is good if
$X_{K}(\underline{x})=1$; it is bad otherwise. According to lemma
\ref{evtR}, we know that for a fixed $\underline{x}$ the variable
$X_{K}(\underline{x})$ is a Bernoulli random variable with parameter
$1-\delta_{K}$, where $\lim_{K\rightarrow \infty} \delta_{K} =0$. This way we
obtain a dependent percolation by edges on the rescaled lattice.

\begin{figure}[h!]
\centering
\input{erreur3.pstex_t}
\caption{Construction of the open path - 2.}
\label{erreur3}
\end{figure}
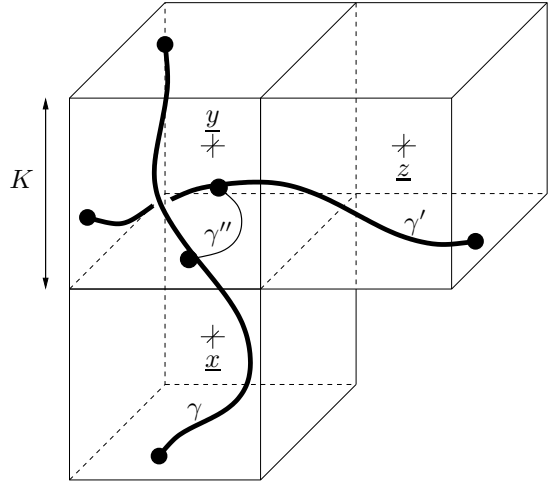

Now if we have a $L^{1}$-connected path of good event-blocks
$(B'_{K}(\underline{x}_{i}), i\in I)$ in the rescaled lattice from the
bottom to the top of a rescaled cylinder $\underline{B}$, we can find an open path from the bottom to the top of the corresponding cylinder $B$
in the initial graph which is completely included in $\cup_{i\in I}
B_{K}(\underline{x}_{i})$. Indeed, take $\underline{x}$, $\underline{y}$ and
$\underline{z}$ three successive elements of the $L^{1}$-connected sequence
$(\underline{x}_{i}, i\in I)$. Suppose (for the recurrence) that we have
already constructed an open path $\gamma$ in $B_{K}(\underline{x}) \cup
B_{K}(\underline{y})$ which join the two opposite faces of $\partial^{in}
(B_{K}(\underline{x}) \cup B_{K}(\underline{y}))$ at distance $2K$ (the
ones perpendicular to the direction of the vector
$\underline{y}-\underline{x}$). We
know that $B_{K}'(\underline{y})$ and $B_{K}'(\underline{z})$ are good
event-blocks, so the events
$U(B_{K}(\underline{y}))$, $U(B_{K}(\underline{z}))$ and
$W(B_{K}(\underline{y})+(\underline{z}-\underline{y}) K/2, K/3)$ occur. On
$U(B_{K}(\underline{y})) \cap U(B_{K}(\underline{z}))$, we know that there exists an open path
$\gamma'_{1}$ (respectively $\gamma'_{2}$) in $B_{K}(\underline{y})$
(respectively $B_{K}(\underline{z})$) that join the two opposite faces of
$\partial^{in} B_{K}(\underline{y})$ (respectively $\partial^{in}
B_{K}(\underline{z})$) perpendicular to the direction of the vector
$\underline{z}-\underline{y}$. Moreover, since
$W(B_{K}(\underline{y})+(\underline{z}-\underline{y}) K/2, K/3)$ occurs, $\gamma'_{1}$ and $\gamma'_{2}$
are connected by an open path $\gamma'_{3}$ in
$B_{K}(\underline{y})+(\underline{z}-\underline{y}) K/2$ because of the
uniqueness of the open cluster of diameter greater than $K/3$ in
$B_{K}(\underline{y})+(\underline{z}-\underline{y}) K/2$ (see figure \ref{erreur2}). So $\gamma' =
\gamma'_{1} \cup \gamma'_{2} \cup \gamma'_{3}$ contains an open path in
$B_{K}(\underline{y}) \cup 
B_{K}(\underline{z})$ which joins the two opposite faces of $\partial^{in}
(B_{K}(\underline{y}) \cup B_{K}(\underline{z}))$ at distance $2K$ (the
ones perpendicular to the direction of the vector
$\underline{z}-\underline{y}$). Finally the event
$W(B_{K}(\underline{y}), K/3)$ occurs so we know that $\gamma$ and
$\gamma'$ are connected by an open path $\gamma''$ in
$B_{K}(\underline{y})$ (see figure \ref{erreur3}). From an event-block to another, we can build the
desired open path from the bottom to the top of the cylinder $B$, and it
lies indeed in $\cup_{i\in I} B_{K}(\underline{x}_{i})$.

\begin{figure}[h!]
\centering
\input{erreur2.pstex_t}
\caption{Construction of the open path - 1.}
\label{erreur2}
\end{figure}
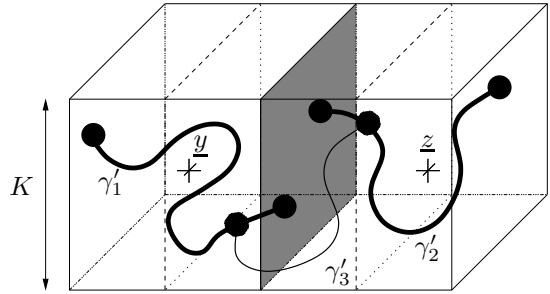

Moreover, if we have $N$ disjoint $L^{1}$-connected paths of good event-blocks
$(B'_{K}(\underline{x}_{i}), i\in I_{j})$, $j=1,...,N$ in the rescaled lattice from the
bottom to the top of a rescaled cylinder $\underline{B}$, we can find N disjoint open paths from the bottom to the top of the corresponding cylinder $B$
in the initial graph, because the sets $I_{j}$ are pairwise disjoint and
so are the sets $\cup_{i\in I_{j}} B_{K}(\underline{x}_{i})$, which contains
the different open paths constructed as previously.


\subsection{Proof of theorem \ref{bernoulli}}
\label{equiv}

As previously we use the notation
$$ \alpha(\varepsilon) \,=\, \mathbb{P} \left[ \phi_{(n,...,n),h(n)} \leq \varepsilon n^{d-1} \right] \,. $$
We define the cylinder
$$ A(n) \,=\, [0,n]^{d-1} \times [0,h(n)] \,,$$
and $\underline{A}_{K}$ is the rescaled cylinder for an integer $K$ which will be chosen soon.

According to the remark at the end of the previous subsection, we know that if there exist $\varepsilon n^{d-1}$ disjoint paths of $L^{1}$-connected good event-blocks from the bottom to the top of the rescaled cylinder $\underline{A}_{K}(n)$, then there exist at least $\varepsilon n^{d-1}$ disjoint open paths from the bottom to the top of $A(n)$. Therefore
$$ \alpha(\varepsilon) \,\leq\, \mathbb{P} \left[ \,\begin{array}{c} there \,\, exist \,\, less \,\, than \,\, \varepsilon n^{d-1} \,\, disjoint \,\, paths \,\, of \,\, good \\  event-blocks \,\, from\,\,the\,\,bottom\,\,to\,\,the\,\,top\,\,of\,\, \underline{A}_{K}(n) \end{array} \, \right] \,.$$
Now the arguments which will be used are very similar to those used in the proof of theorem \ref{bernoulli1}. The main difference is that the random variables $(X_{K}(\underline{x}), \underline{x}\in\mathbb{Z}^{d})$ are not independent.

We work for the rest of the proof in the rescaled lattice. The notion of cut can be adapted easily in the model of site percolation and the max-flow min-cut theorem remains valid in this model. Thanks to the max-flow min-cut theorem applied in the rescaled lattice we obtain that
$$ \alpha(\varepsilon) \,\leq\, \mathbb{P} \left[\, there \,\, exists \,\, a \,\, (\underline{F_{0}}_{K},\underline{F_{h(n)}}_{K})-cut \,\, \underline{E} \,\, satisfying \,\, \underline{V}(\underline{E})\leq \varepsilon n^{d-1}\,  \right] \,,$$
where here
$$ \underline{V}(\underline{E}) \,=\, \sum_{\underline{x}\in\underline{E}} X_{K}(\underline{x}) \,. $$
Note that such a $(\underline{F_{0}}_{K},\underline{F_{h(n)}}_{K})$-cut contains at least $u=\lfloor (n/K)^{d-1} \rfloor$ vertices. As previously, we obtain
$$ \alpha(\varepsilon) \leq\, \frac{h(n)}{K} \sum_{\underline{A}} \mathbb{P} \left[ \sum_{\underline{x}\in\underline{A}} X_{K}(\underline{x}) \leq \varepsilon n^{d-1} \right] \, ,$$
where the sum is over the $L^{1}$-connected sets $\underline{A}$ of $u$ vertices containing a fixed vertex $\underline{x}_{0}$ of $\mathbb{Z}^{d}$.

To deal with the variables $(X_{K}(\underline{x}), \underline{x}\in\mathbb{Z}^{d})$ we introduce an equivalence relation on $\mathbb{Z}^{d}$: $x\sim y$ if and only if $3$ divides all the coordinates of $x-y$. There exist $3^{d}$ equivalence classes $V_{1},...,V_{3^{d}}$ in $\mathbb{Z}^{d}$. For a set of vertices $\underline{E}$, we define
$$ \underline{E}^{l} \,=\, \underline{E} \cap V_{l} \,. $$
Now the variables $(X_{K}(\underline{x}), \underline{x}\in V_{l})$ are independent for a fixed $l \in \{1,...,3^{d}\}$, so we want to consider only sums of variables indexed by vertices in the same equivalence class. For that purpose, we remark that if $\sum_{\underline{x}\in \underline{A}} X_{K}(\underline{x}) \leq \varepsilon n^{d-1}$ for some set $\underline{A}$ of $u$ vertices and for some $\varepsilon \leq 1/K^{d-1}$, then $\underline{A}$ contains at least $u - \lfloor  \varepsilon n^{d-1} \rfloor$ bad event-blocks which are included in the subsets $\underline{A}^{1},...,\underline{A}^{3^{d}}$. Thus there exists $l\in \{1,...,3^{d}\}$ such that $\underline{A}^{l}$ contains at least $(u - \lfloor \varepsilon n^{d-1} \rfloor )/3^{d}$ bad event-blocks, so $\sum_{\underline{x}\in \underline{A}^{l}} X_{K}(\underline{x}) \leq |\underline{A}^{l}| -(u - \lfloor \varepsilon n^{d-1} \rfloor )/3^{d} $. This remark leads to
$$ \mathbb{P} \left[ \sum_{\underline{x}\in\underline{A}} X_{K}(\underline{x}) \leq \varepsilon n^{d-1} \right] \,\leq\, \sum_{l=1}^{3^{d}} \mathbb{P} \left[ \sum_{\underline{x}\in \underline{A}^{l}} X_{K}(\underline{x}) \leq |\underline{A}^{l}| - \frac{1}{3^{d}} (u-\varepsilon n^{d-1}) \right] \,, $$
where $| \underline{A}^{l}|$ is the cardinal of $\underline{A}^{l}$. Now, thanks again to the bound on the number of possible sets $\underline{A}$ and to the exponential Chebyshev inequality, we obtain as in the proof of theorem \ref{bernoulli1}
\begin{align*}
\forall \lambda >0 \qquad \alpha (\varepsilon) & \,\leq\, \frac{h(n)}{K} \sum_{\underline{A}} \sum_{l=1}^{3^{d}} \exp \left( \lambda \left[ |\underline{A}^{l}| - \frac{1}{3^{d}} (u-\varepsilon n^{d-1}) \right] \right) \mathbb{E} \left[ e^{-\lambda X_{K}(\underline{x}_{0})} \right]^{|\underline{A}^{l}|} \\
& \,\leq\, \frac{h(n)}{K} \sum_{\underline{A}} \sum_{l=1}^{3^{d}} \exp \left( \lambda |\underline{A}^{l}| - \lambda \frac{u-\varepsilon n^{d-1}}{3^{d}} + |\underline{A}^{l}| \ln \left[ e^{-\lambda}(1+ \delta_{K} (e^{\lambda}-1))  \right]  \right) \\
& \,\leq\, \frac{3^{d}h(n)}{K} \exp \left( -u \left[-\ln c' +\frac{\lambda}{3^{d}} (1- \frac{\varepsilon n^{d-1}}{u}) - \ln (1+\delta_{K} (e^{\lambda}-1))  \right] \right) \,, \\
\end{align*}
because $|\underline{A}^{l}|\leq u$. Now we choose first $\lambda $ such that
$$ \frac{\lambda}{2\times 3^{d}} \, \geq \, 3 \ln c' \,, $$
and then $K$ large enough (depending on $d$ and $p$), so $\delta_{K}$ small enough, to have
$$ \ln \left( 1+ \delta_{K} (e^{\lambda} -1) \right)\, \leq \,\ln c' \,. $$
We obtain
$$ \forall \varepsilon \leq \frac{1}{2K^{d-1}} \qquad \alpha(\varepsilon)\, \leq \, \frac{3^{d}h(n)}{K} e^{-u \ln c'} \,, $$
and this ends the proof.


\section[Proof of theorem 1]{Proof of theorem \ref{casgen}}

We consider finally the general case of the first passage percolation model, with the condition
$$ F(0) \,<\, 1-p_{c} \,.$$
The distribution function $F$ is right continuous, so there exists a positive $\eta$ such that
$$p' \,=\, \mathbb{P} [t > \eta] \,=\, 1 -F(\eta) \,>\, p_{c} \,.$$
Now we consider a new family of random variables on $\mathbb{E}^{d}$ defined as
$$ t'(e) \,=\, \left\{ \begin{array}{ll} 1 & \,\, if \,\, t(e)>\eta \\ 0 & \,\, otherwise \end{array}  \right. $$
The family $(t'(e), e\in\mathbb{E}^{d})$ defines an independent Bernoulli percolation of parameter $p'$ on the lattice. We consider the rescaled lattice, and we say that an event-block is good if it is good for this Bernoulli percolation according to the definition given in the previous section. We remark that the existence of a path of good event-blocks in the rescaled lattice implies the existence of a path of edges with a capacity greater than $\eta$ in the initial graph. Therefore
\begin{align*}
\alpha (\varepsilon) & \,=\, \mathbb{P} \left[ \phi_{(n,...,n),h(n)} \leq \varepsilon n^{d-1} \right] \\
& \,\leq\, \mathbb{P} \left[ \,\begin{array}{c} there \,\, exist \,\, less \,\, than \,\,\frac{ \varepsilon}{\eta} n^{d-1} \,\, disjoint \,\, paths \,\, of \,\, good \\ event-blocks \,\, from\,\,the\,\,bottom\,\,to\,\,the\,\,top\,\, of\,\, \underline{A}_{K}(n)\end{array} \,  \right] \,. \\
\end{align*}
We proceed as in the proof of theorem \ref{bernoulli} to obtain the desired estimate.
\\

{\bf Acknowledgement:} The author would like to warmly thank Rapha\"el Cerf for his guidance, his help and his kindness. The author is also very
grateful to Rapha\"el Rossignol for his numerous valuable comments, one of
which revealed a mistake in the construction of the good event-blocks.


\end{document}

%% file: erreur2.pstex_t
\begin{picture}(0,0)%
\includegraphics{erreur2.pstex}%
\end{picture}%
\setlength{\unitlength}{1973sp}%
\begingroup\makeatletter\ifx\SetFigFont\undefined%
\gdef\SetFigFont#1#2#3#4#5{%
  \reset@font\fontsize{#1}{#2pt}%
  \fontfamily{#3}\fontseries{#4}\fontshape{#5}%
  \selectfont}%
\fi\endgroup%
\begin{picture}(6861,3624)(1254,-5173)
\put(6601,-3361){\makebox(0,0)[b]{\smash{{\SetFigFont{10}{12.0}{\rmdefault}{\mddefault}{\updefault}{\color[rgb]{0,0,0}$\underline{z}$}%
}}}}
\put(3751,-3361){\makebox(0,0)[b]{\smash{{\SetFigFont{10}{12.0}{\rmdefault}{\mddefault}{\updefault}{\color[rgb]{0,0,0}$\underline{y}$}%
}}}}
\put(2626,-3886){\makebox(0,0)[b]{\smash{{\SetFigFont{10}{12.0}{\rmdefault}{\mddefault}{\updefault}{\color[rgb]{0,0,0}$\gamma'_{1}$}%
}}}}
\put(5476,-5011){\makebox(0,0)[b]{\smash{{\SetFigFont{10}{12.0}{\rmdefault}{\mddefault}{\updefault}{\color[rgb]{0,0,0}$\gamma'_{3}$}%
}}}}
\put(6601,-4711){\makebox(0,0)[b]{\smash{{\SetFigFont{10}{12.0}{\rmdefault}{\mddefault}{\updefault}{\color[rgb]{0,0,0}$\gamma'_{2}$}%
}}}}
\put(1501,-3961){\makebox(0,0)[b]{\smash{{\SetFigFont{10}{12.0}{\rmdefault}{\mddefault}{\updefault}{\color[rgb]{0,0,0}$K$}%
}}}}
\end{picture}%

%% file: erreur3.pstex_t
\begin{picture}(0,0)%
\includegraphics{erreur3.pstex}%
\end{picture}%
\setlength{\unitlength}{1973sp}%
\begingroup\makeatletter\ifx\SetFigFont\undefined%
\gdef\SetFigFont#1#2#3#4#5{%
  \reset@font\fontsize{#1}{#2pt}%
  \fontfamily{#3}\fontseries{#4}\fontshape{#5}%
  \selectfont}%
\fi\endgroup%
\begin{picture}(6859,6024)(1254,-6973)
\put(3901,-5536){\makebox(0,0)[b]{\smash{{\SetFigFont{10}{12.0}{\rmdefault}{\mddefault}{\updefault}{\color[rgb]{0,0,0}$\underline{x}$}%
}}}}
\put(3901,-2461){\makebox(0,0)[b]{\smash{{\SetFigFont{10}{12.0}{\rmdefault}{\mddefault}{\updefault}{\color[rgb]{0,0,0}$\underline{y}$}%
}}}}
\put(6301,-3136){\makebox(0,0)[b]{\smash{{\SetFigFont{10}{12.0}{\rmdefault}{\mddefault}{\updefault}{\color[rgb]{0,0,0}$\underline{z}$}%
}}}}
\put(1501,-3286){\makebox(0,0)[b]{\smash{{\SetFigFont{10}{12.0}{\rmdefault}{\mddefault}{\updefault}{\color[rgb]{0,0,0}$K$}%
}}}}
\put(3976,-3961){\makebox(0,0)[b]{\smash{{\SetFigFont{10}{12.0}{\rmdefault}{\mddefault}{\updefault}{\color[rgb]{0,0,0}$\gamma''$}%
}}}}
\put(3676,-6136){\makebox(0,0)[b]{\smash{{\SetFigFont{10}{12.0}{\rmdefault}{\mddefault}{\updefault}{\color[rgb]{0,0,0}$\gamma$}%
}}}}
\put(6451,-3811){\makebox(0,0)[b]{\smash{{\SetFigFont{10}{12.0}{\rmdefault}{\mddefault}{\updefault}{\color[rgb]{0,0,0}$\gamma'$}%
}}}}
\end{picture}%